\documentclass[reqno,12pt,a4paper]{amsart}

\usepackage{mathrsfs}
\usepackage{amssymb}
\usepackage{amsfonts}
\usepackage{latexsym}
\usepackage{amsthm}

\usepackage{graphicx}
\def\lb{\label}

\newcommand{\er}[1]{\textrm{(\ref{#1})}}

\begin{document}

%%%%%%%%%% Some definitions %%%%%%%%%%

%%%%%%%% Equations, theorems %%%%%%%%%
\renewcommand{\theequation}{\arabic{section}.\arabic{equation}}
\theoremstyle{plain}
\newtheorem{theorem}{\bf Theorem}[section]
\newtheorem{lemma}[theorem]{\bf Lemma}
\newtheorem{corollary}[theorem]{\bf Corollary}
\newtheorem{proposition}[theorem]{\bf Proposition}
\newtheorem{definition}[theorem]{\bf Definition}
\newtheorem{remark}[theorem]{\it Remark}
%\theoremstyle{remark}
%\newtheorem{remark}[theorem]{\bf Remark}

%%%%% Alphabet %%%%%
\def\a{\alpha}  \def\cA{{\mathcal A}}     \def\bA{{\bf A}}  \def\mA{{\mathscr A}}
\def\b{\beta}   \def\cB{{\mathcal B}}     \def\bB{{\bf B}}  \def\mB{{\mathscr B}}
\def\g{\gamma}  \def\cC{{\mathcal C}}     \def\bC{{\bf C}}  \def\mC{{\mathscr C}}
\def\G{\Gamma}  \def\cD{{\mathcal D}}     \def\bD{{\bf D}}  \def\mD{{\mathscr D}}
\def\d{\delta}  \def\cE{{\mathcal E}}     \def\bE{{\bf E}}  \def\mE{{\mathscr E}}
\def\D{\Delta}  \def\cF{{\mathcal F}}     \def\bF{{\bf F}}  \def\mF{{\mathscr F}}
\def\c{\chi}    \def\cG{{\mathcal G}}     \def\bG{{\bf G}}  \def\mG{{\mathscr G}}
\def\z{\zeta}   \def\cH{{\mathcal H}}     \def\bH{{\bf H}}  \def\mH{{\mathscr H}}
\def\e{\eta}    \def\cI{{\mathcal I}}     \def\bI{{\bf I}}  \def\mI{{\mathscr I}}
\def\p{\psi}    \def\cJ{{\mathcal J}}     \def\bJ{{\bf J}}  \def\mJ{{\mathscr J}}
\def\vT{\Theta} \def\cK{{\mathcal K}}     \def\bK{{\bf K}}  \def\mK{{\mathscr K}}
\def\k{\kappa}  \def\cL{{\mathcal L}}     \def\bL{{\bf L}}  \def\mL{{\mathscr L}}
\def\l{\lambda} \def\cM{{\mathcal M}}     \def\bM{{\bf M}}  \def\mM{{\mathscr M}}
\def\L{\Lambda} \def\cN{{\mathcal N}}     \def\bN{{\bf N}}  \def\mN{{\mathscr N}}
\def\m{\mu}     \def\cO{{\mathcal O}}     \def\bO{{\bf O}}  \def\mO{{\mathscr O}}
\def\n{\nu}     \def\cP{{\mathcal P}}     \def\bP{{\bf P}}  \def\mP{{\mathscr P}}
\def\r{\rho}    \def\cQ{{\mathcal Q}}     \def\bQ{{\bf Q}}  \def\mQ{{\mathscr Q}}
\def\s{\sigma}  \def\cR{{\mathcal R}}     \def\bR{{\bf R}}  \def\mR{{\mathscr R}}
                \def\cS{{\mathcal S}}     \def\bS{{\bf S}}  \def\mS{{\mathscr S}}
\def\t{\tau}    \def\cT{{\mathcal T}}     \def\bT{{\bf T}}  \def\mT{{\mathscr T}}
\def\f{\phi}    \def\cU{{\mathcal U}}     \def\bU{{\bf U}}  \def\mU{{\mathscr U}}
\def\F{\Phi}    \def\cV{{\mathcal V}}     \def\bV{{\bf V}}  \def\mV{{\mathscr V}}
\def\P{\Psi}    \def\cW{{\mathcal W}}     \def\bW{{\bf W}}  \def\mW{{\mathscr W}}
\def\o{\omega}  \def\cX{{\mathcal X}}     \def\bX{{\bf X}}  \def\mX{{\mathscr X}}
\def\x{\xi}     \def\cY{{\mathcal Y}}     \def\bY{{\bf Y}}  \def\mY{{\mathscr Y}}
\def\X{\Xi}     \def\cZ{{\mathcal Z}}     \def\bZ{{\bf Z}}  \def\mZ{{\mathscr Z}}
\def\O{\Omega}

\newcommand{\gA}{\mathfrak{A}}          \newcommand{\ga}{\mathfrak{a}}
\newcommand{\gB}{\mathfrak{B}}          \newcommand{\gb}{\mathfrak{b}}
\newcommand{\gC}{\mathfrak{C}}          \newcommand{\gc}{\mathfrak{c}}
\newcommand{\gD}{\mathfrak{D}}          \newcommand{\gd}{\mathfrak{d}}
\newcommand{\gE}{\mathfrak{E}}
\newcommand{\gF}{\mathfrak{F}}           \newcommand{\gf}{\mathfrak{f}}
\newcommand{\gG}{\mathfrak{G}}           %\newcommand{\gg}{\mathfrak{g}}
\newcommand{\gH}{\mathfrak{H}}           \newcommand{\gh}{\mathfrak{h}}
\newcommand{\gI}{\mathfrak{I}}           \newcommand{\gi}{\mathfrak{i}}
\newcommand{\gJ}{\mathfrak{J}}           \newcommand{\gj}{\mathfrak{j}}
\newcommand{\gK}{\mathfrak{K}}            \newcommand{\gk}{\mathfrak{k}}
\newcommand{\gL}{\mathfrak{L}}            \newcommand{\gl}{\mathfrak{l}}
\newcommand{\gM}{\mathfrak{M}}            \newcommand{\gm}{\mathfrak{m}}
\newcommand{\gN}{\mathfrak{N}}            \newcommand{\gn}{\mathfrak{n}}
\newcommand{\gO}{\mathfrak{O}}
\newcommand{\gP}{\mathfrak{P}}             \newcommand{\gp}{\mathfrak{p}}
\newcommand{\gQ}{\mathfrak{Q}}             \newcommand{\gq}{\mathfrak{q}}
\newcommand{\gR}{\mathfrak{R}}             \newcommand{\gr}{\mathfrak{r}}
\newcommand{\gS}{\mathfrak{S}}              \newcommand{\gs}{\mathfrak{s}}
\newcommand{\gT}{\mathfrak{T}}             \newcommand{\gt}{\mathfrak{t}}
\newcommand{\gU}{\mathfrak{U}}             \newcommand{\gu}{\mathfrak{u}}
\newcommand{\gV}{\mathfrak{V}}             \newcommand{\gv}{\mathfrak{v}}
\newcommand{\gW}{\mathfrak{W}}             \newcommand{\gw}{\mathfrak{w}}
\newcommand{\gX}{\mathfrak{X}}               \newcommand{\gx}{\mathfrak{x}}
\newcommand{\gY}{\mathfrak{Y}}              \newcommand{\gy}{\mathfrak{y}}
\newcommand{\gZ}{\mathfrak{Z}}             \newcommand{\gz}{\mathfrak{z}}

\def\ve{\varepsilon} \def\vt{\vartheta} \def\vp{\varphi}  \def\vk{\varkappa}

\def\Z{{\mathbb Z}} \def\R{{\mathbb R}} \def\C{{\mathbb C}}  \def\K{{\mathbb K}}
\def\T{{\mathbb T}} \def\N{{\mathbb N}} \def\dD{{\mathbb D}} \def\S{{\mathbb S}}
\def\B{{\mathbb B}}

%%%%% Arrows %%%%%

\def\la{\leftarrow}              \def\ra{\rightarrow}            \def\Ra{\Rightarrow}
\def\ua{\uparrow}                \def\da{\downarrow}
\def\lra{\leftrightarrow}        \def\Lra{\Leftrightarrow}

%%%%% Typography %%%%%

\def\lt{\biggl}                  \def\rt{\biggr}
\def\ol{\overline}               \def\wt{\widetilde}
\def\no{\noindent}

%%%%% Math signs %%%%%

\let\ge\geqslant                 \let\le\leqslant
\def\lan{\langle}                \def\ran{\rangle}
\def\/{\over}                    \def\iy{\infty}
\def\sm{\setminus}               \def\es{\emptyset}
\def\ss{\subset}                 \def\ts{\times}
\def\pa{\partial}                \def\os{\oplus}
\def\om{\ominus}                 \def\ev{\equiv}
\def\iint{\int\!\!\!\int}        \def\iintt{\mathop{\int\!\!\int\!\!\dots\!\!\int}\limits}
\def\el2{\ell^{\,2}}             \def\1{1\!\!1}
\def\sh{\sharp}
\def\wh{\widehat}
\def\bs{\backslash}
%%%%% Math operations %%%%%

\def\sh{\mathop{\mathrm{sh}}\nolimits}
\def\Area{\mathop{\mathrm{Area}}\nolimits}
\def\arg{\mathop{\mathrm{arg}}\nolimits}
\def\const{\mathop{\mathrm{const}}\nolimits}
\def\det{\mathop{\mathrm{det}}\nolimits}
\def\diag{\mathop{\mathrm{diag}}\nolimits}
\def\diam{\mathop{\mathrm{diam}}\nolimits}
\def\dim{\mathop{\mathrm{dim}}\nolimits}
\def\dist{\mathop{\mathrm{dist}}\nolimits}
\def\Im{\mathop{\mathrm{Im}}\nolimits}
\def\Iso{\mathop{\mathrm{Iso}}\nolimits}
\def\Ker{\mathop{\mathrm{Ker}}\nolimits}
\def\Lip{\mathop{\mathrm{Lip}}\nolimits}
\def\rank{\mathop{\mathrm{rank}}\limits}
\def\Ran{\mathop{\mathrm{Ran}}\nolimits}
\def\Re{\mathop{\mathrm{Re}}\nolimits}
\def\Res{\mathop{\mathrm{Res}}\nolimits}
\def\res{\mathop{\mathrm{res}}\limits}
\def\sign{\mathop{\mathrm{sign}}\nolimits}
\def\span{\mathop{\mathrm{span}}\nolimits}
\def\supp{\mathop{\mathrm{supp}}\nolimits}
\def\Tr{\mathop{\mathrm{Tr}}\nolimits}
\def\BBox{\hspace{1mm}\vrule height6pt width5.5pt depth0pt \hspace{6pt}}
\def\as{\text{as}}
\def\all{\text{all}}
\def\where{\text{where}}
\def\Dom{\mathop{\mathrm{Dom}}\nolimits}
\def\ch{\mathop{\mathrm{ch}}\nolimits}
\def\sh{\mathop{\mathrm{sh}}\nolimits}

%%%%%%%%%%%%% specialities %%%%%%%%%%%%%%

\newcommand\nh[2]{\widehat{#1}\vphantom{#1}^{(#2)}}
%{{\mathop{#1}\limits^\wedge}\vphantom{#1}^{(#2)}}
\def\dia{\diamond}

\def\Oplus{\bigoplus\nolimits}

%%%%%%%%%%% End of definitions %%%%%%%%%%

%%%%% OLD OLD OLD

\def\qqq{\qquad}
\def\qq{\quad}
\let\ge\geqslant
\let\le\leqslant
\let\geq\geqslant
\let\leq\leqslant
\newcommand{\ca}{\begin{cases}}
\newcommand{\ac}{\end{cases}}
\newcommand{\ma}{\begin{pmatrix}}
\newcommand{\am}{\end{pmatrix}}
\renewcommand{\[}{\begin{equation}}
\renewcommand{\]}{\end{equation}}
\def\eq{\begin{equation}}
\def\qe{\end{equation}}
\def\[{\begin{equation}}
\def\bu{\bullet}
\newcommand{\fr}{\frac}
\newcommand{\tf}{\tfrac}

\title[{Schr\"odinger operators  periodic in octants}]
{Schr\"odinger operators  periodic in octants}

\date{\today}
\author[Evgeny Korotyaev]{Evgeny Korotyaev}
\address{Saint-Petersburg State University,
Universitetskaya nab. 7/9, St. Petersburg, 199034, Russia,
\ korotyaev@gmail.com, \
e.korotyaev@spbu.ru,}

\author[Jacob  Schach  Moller]{Jacob  Schach  Moller}
\address{Institut for Matematiske Fag, Aarhus Universitet,
Ny Munkegade 8000 Aarhus C, Denmark, \ jacob@math.au.dk}


\begin{abstract}
We consider Schr\"odinger operators with  periodic potentials in the
positive  quadrant for dim $>1$ with Dirichlet boundary condition.
We show that for any integer $N$ and any interval $I$ there exists a
periodic potential such that the Schr\"odinger operator has $N$
eigenvalues counted with the multiplicity on this interval and there
is no other spectrum on the interval. Furthermore, to the right and
to the left of it there is a essential spectrum.
 Moreover, we prove
similar results for Schr\"odinger operators for other domains. The
proof is based on the inverse spectral theory for Hill operators on
the real line.

\end{abstract}

\maketitle

\vskip 0.25cm

\section {\lb{Sec1}Introduction and main results}
\setcounter{equation}{0}

\subsection{Introduction.}

We consider  Schr\"odinger operators $H$ on the domain $D$
 given by
\[
\lb{d1}
\begin{aligned}
& H=-\D_x^+-\D_y+V(x,y),\\
&  (x,y)\in D=\R_+^{d_1}\ts \R^{d_2},\qq d_1+d_2=d\ge 2,\qq
d_1,d_2\ge
  0.
\end{aligned}
\]
%where $d_1+d_2=d\ge 2$ and $d_1,d_2\ge 0$.
Here the operator
$\D_x^+$ is the Laplacian in the octant  $ \R_+^{d_1}$ with the
Dirichlet boundary conditions on the boundary $\pa \R_+^{d_1}$ and
the operator $\D_y$ is the Laplacian in  the space $ \R^{d_2}$. We
assume that the potential $V$ belongs to $ L_{real}^\iy(D)$ and is
octant periodic, see Condition V.

In order to define octant periodic potentials  we need additional
definitions. Let $\o=(\o_j)_1^m$ be a sequence of $+$ or $-$  and
the set all these sequences we denote by $\O_m$. For any $\o\in
\O_m$ we define the octants $\cR_\o\ss \R^{m}$ by
$$
\cR_\o=\R_{\o_1}\ts \R_{\o_2}\ts ....\ts
\R_{\o_m},\qqq\o=(\o_j)_1^m\in \O_m.
$$
In particular, if  $d=3$ and $\o=(+,+,+)$, then we have the positive
octant $\cR_\o=\R_+^3$. Note that three axial planes $(x_1=0, x_2=0,
x_3=0)$ divide space $\R^3$ into eight octants, each with a
coordinate signs from $(-,-,-)$ to $(+,+,+)$.

{\bf Definition V.} {\it  A potential $V(z), z=(x,y)\in
\R_+^{d_1}\ts \R^{d_2}$ is  called octant periodic if it
  has the decomposition
\[
\lb{VO}
\begin{aligned}
  V(x,y)=\sum_{\o\in \O_{d_2}} V_\o(x,y)\c_\o(y),\qqq
  %& V_\o(x+\gp_{\o,j}{\bf e}_j)=V_\o(x),\qq \forall \ x\in \R^d,\qq
%j=1,...,d,
\end{aligned}
\]
where $\c_\o$ is the characteristic function of the octant $\cR_\o$
and the function $V_\o(z), z=(x,y)$ is periodic in $\R^d$ and
satisfies
\[
\lb{V}
\begin{aligned}
 V(z+\gp_j{\bf e}_j)=V(z),\qq \forall \ z\in \R^d,\qq
j=1,...,d,
\end{aligned}
\]
for some constants $\gp_j>0, j=1,2,..,d$, where  and ${\bf
e}_1=(1,0,0,..), {\bf e}_2=(0,1,0,0,..)$ is the standard basis in
$\R^d$.}

For each $\o$ we define  Schr\"odinger operators $H_\o$ with
periodic potentials $V_\o$ on $\R^d$  by
\[
\lb{dHo}
\begin{aligned}
 H_\o=-\D_z+V_\o(z).
\end{aligned}
\]
It is well known that the spectrum of each operator $H_\o$ is
absolutely continuous and  is an union of an unbounded interval and
a finite number of non-degenerated bounded intervals. In the next
theorem we show the existence of eigenvalues of $H$ with some octant
periodic potentials.

\begin{theorem}
\lb{T1} i) Let $H=-\D_x^+-\D_y+V(x,y)$, where $(x,y)\in
D=\R_+^{d_1}\ts \R^{d_2}, d_1,d_2\ge 0, d_1+d_2=d\ge 2$ and the
potential $V\in L^\iy(D)$ is octant periodic. Then
\[
\lb{spH} \bigcup_{\o\in \O_{d_2}}\s(H_\o) \subseteq \s_{ess}(H).
\]
ii)  Let $I=(a,b)\ss \R$ be a finite open interval. Then for any
integer $N\ge 0$ there exists an octant periodic potential $V\in
L^\iy(D)$
 such that $H=-\D_x^+-\D_y+V(z)$ on $D$ has $N$
eigenvalues counted with multiplicity on $I$. Moreover, the interval
$I$ does not contain the essential spectrum, to the right and to the
left of it there is a essential spectrum.
%each of the intervals $[a_-,a]$ and $[b,b_+]$ contains the essential
%spectrum of $H$ for some $a_-,b_+$.
\end{theorem}

{\bf Remark.} 1) The relation \er{spH} is standard. Its proof
follows from the Floquet theory for periodic operators.

2)  In the simple case we have Schr\"odinger operators
$H=-\D_x^++V(x)$ on the quadrant $\R_+^2$, where the operator
$\D_x^+$ is the Laplacian in  the domain $ \R_+^2$ with the
Dirichlet boundary conditions on the boundary $\pa \R_+^2$. We
assume that the potential $V$ is octant periodic. This theorem
describe the spectrum of $H$ on the positive quadrant.

3) The proof of ii) is based on the inverse spectral theory from
\cite{K99}.

Now we discuss an existence of eigenvalue  of the operator $H$
acting on $\R^d$ below the essential spectrum.

\begin{theorem}
\lb{T2} There exists a octant periodic real potential $V$ such that
the operator $H=-\D_y+V(y)$ on $\R^d$ has exactly one simple
eigenvalue on some interval $(-\iy, E)$ and there are no other
spectrum on it.

\end{theorem}

{\bf Remark.} There is a problem to show that there is only one
eigenvalue below the continuous spectrum.

\subsection{Historical review}

Davies and Simon \cite{DS78} study the Schr\"odinger operators with
potential, which periodic in half-space.  The model one dimensional
case with periodic potentials on the half-line is considered in
\cite{K00}, \cite{K05}. In series of papers
\cite{HK11},\cite{HK11x},\cite{HKSV15}, Hempel and Kohlmann with
co-authers consider the different type of dislocation problem in
solid state physics. The surface density of states is investigated
in \cite{JL01}, \cite{KS01}. Scattering on the periodic boundary is
considered in \cite{F06}, \cite{JL03}.

%\bigskip

In the simple case $d=1$  we need to describe the spectrum of two
operators:

$\bu$ In Section 2 we consider the Schr\"odinger operator $h$  on
the space $L^2(\R_+ )$ given by
$$
hy=-y''+vy, \qqq
$$
with the boundary condition $y(0)=0$. Here the potential $v$ is real
1-periodic and $ v\in L^1(0,1)$. In our proof results from
\cite{K99}, \cite{K06} about the inverse problem for $h$ from  is
crucial and is presented in Theorem \ref{Ti}.
%We recall the result from \cite{K05}.

$\bu$  We also consider the so-called half-solid operator $T_\t,
\t\in \R$ acting on $L^2(\R )$ and given by
\[
\lb{qsx}
 \ T_\t=-{d^2\/dx^2}+q_\t(x),\qqq q_\t(x)=\ca \t & if \ x<0\\
           v(x) & if \ x>0\ac,\
\]
where the potential $v\in L_{real}^2(\T)$. By the physical point of
view $v$ is the potential of a crystal and the constant $\t$ is the
potential of a vacuum. Roughly speaking the potentia $q_\t$ is the
octant periodic and has the form $q_\t=V_+\c_++V_-\c_-$, where
$V_+=v $ and $V_-=\t$ and $\c_\pm$ are the characteristic functions
of the half-line $\R_\pm$. In our proof results from \cite{K05}
about spectral properties of $h$ from  is crucial and is presented
in Theorem \ref{Ti}. Moreover, we need additional results formulated
in Lemmas \ref{TD1} and \ref{TD2}.

\section {Periodic Schr\"odinger operators on the half-line}
\setcounter{equation}{0}

\subsection{Preleminary}

We consider the Schr\"odinger operator $h$ acting on the space
$L^2(\R_+ )$ and  given by
$$
hy=-y''+vy, \qqq
$$
with the boundary condition $y(0)=0$. Here the potential $v$ is real
1-periodic and satisfies $ v\in L^1(0,1)$.
 The spectrum of $h$ consists of an
absolutely continuous part $\s_{ac}(h)$ (the union of the bands
$\s_n, n\ge 0$ separated by gaps $\g_n$) plus at most one eigenvalue
in each non-empty gap $\g_n$, $n\in\N$, \cite{KS12}, \cite{Z69}.
Here the bands $\s_n$ and gaps $\g_n$ are given by (see Fig. 1)
\[
\lb{sh}\s_{ac}(h)= \bigcup\limits_{n\ge 0} \s_n,\qq
\s_n=[\l^+_{n},\l^-_{n+1}],\ \ \qq \g_{n}=(\l^-_{n},\l^+_n),\qq
n\in\N,\qq
\]
We also set $\g_0=(-\iy,\l_0^+)$ and $\l_0^+=0$. Here  $\l_n^\pm$
satisfy
\[
\lb{00x}
\begin{aligned}
\l_0^+< \l^-_1 \le \l^+_1\dots\le \l^+_{n-1}< \l^-_n \le
\l^+_{n}<\dots,\\
\l_{n}^\pm=(\pi n)^2+O(1)\qqq \as \qqq n\to \iy.
\end{aligned}
\]
The bands $\s_n, n\ge 0$ satisfy (see e.g. \cite{M81} or \cite{K97})
\[
\lb{eb} |\s_n|=\l^-_{n+1}-\l^+_{n}\le \pi^2(2n+1),\qqq \forall \
n\ge 0.
\]
%It is known that there are infinitely many non-degenerate gaps, i.e.
%$\l_n^- < \l_n^+$, unless $p$ is arbitrarily often differentiable
%%\cite{Ho}, and all gaps are non-degenerate generically \cite{MO},\cite{Si}.
%Without loss of generality, we may assume $\l_0^+=0$.
The sequence \er{00x} is the spectrum of the equation
\[
\lb{1} -y''+v(x)y=\l y, \qqq \l\in \C ,
\]
with the 2-periodic condition $y(x+2)=y(x)$ $(x\in \R)$. If a gap
degenerates, $\g_n=\es $ for some $n\ge 1$, then the corresponding
bands $\s_{n} $ and $\s_{n+1}$ touch. This happens when
$\l_n^-=\l_n^+$; this number is then a double eigenvalue of the
2-periodic problem \er{1}. The lowest eigenvalue $\l_0^+=0$ is
always simple and has a 1-periodic eigenfunction. Generally, the
eigenfunctions corresponding to eigenvalues $\l_{2n}^{\pm}$ are
1-periodic, those for $\l_{2n+1}^{\pm}$ are 1-anti-periodic in the
sense that $y(x+1)=-y(x)$ $(x\in\R)$.

\begin{figure}
\tiny \unitlength=1.00mm \special{em:linewidth 0.4pt}
\linethickness{0.4pt}
\begin{picture}(108.67,33.67)
%coordinate lines
\put(41.00,17.33){\line(1,0){67.67}}
\put(44.33,9.00){\line(0,1){24.67}}
\put(108.33,14.00){\makebox(0,0)[cc]{$\Re\l$}}
\put(41.66,33.67){\makebox(0,0)[cc]{$\Im\l$}}
%\put(42.00,14.33){\makebox(0,0)[cc]{$0$}}
%spectrum
\put(44.33,17.33){\linethickness{4.0pt}\line(1,0){11.33}}
%\put(48.66,17.33){\linethickness{4.0pt}\line(1,0){11.33}}
\put(66.66,17.33){\linethickness{4.0pt}\line(1,0){11.67}}
\put(82.00,17.33){\linethickness{4.0pt}\line(1,0){12.00}}
\put(95.66,17.33){\linethickness{4.0pt}\line(1,0){11.00}}
%endpoints of gaps
\put(46.66,20.00){\makebox(0,0)[cc]{$\l_0^+$}}
%\put(48.66,20.00){\makebox(0,0)[cc]{$E_0^+$}}
\put(56.66,20.33){\makebox(0,0)[cc]{$\l_1^-$}}
%\put(59.66,20.33){\makebox(0,0)[cc]{$E_1^-$}}
\put(68.66,20.33){\makebox(0,0)[cc]{$\l_1^+$}}
%\put(66.66,20.33){\makebox(0,0)[cc]{$E_1^+$}}
\put(78.33,20.33){\makebox(0,0)[cc]{$\l_2^-$}}
\put(84.33,20.33){\makebox(0,0)[cc]{$\l_2^+$}}
%\put(82.33,20.33){\makebox(0,0)[cc]{$E_2^+$}}
\put(93.00,20.33){\makebox(0,0)[cc]{$\l_3^-$}}
\put(98.66,20.33){\makebox(0,0)[cc]{$\l_3^+$}}
%\put(96.66,20.33){\makebox(0,0)[cc]{$E_3^+$}}
\put(106.33,20.33){\makebox(0,0)[cc]{$\l_4^-$}}
\end{picture}
\caption{The cut domain $\C\sm \cup \gS_n$ and the cuts (bands)
$\gS_n=[E^+_{n-1},E^-_n], n\ge 1$} \lb{sS}
\end{figure}

Introduce the two canonical fundamental solutions $\vt(x,\l)$,
$\vp(x,\l)$  of the equation \er{1}, satisfying the initial
conditions $ \vp'(0,\l)=\vt(0,\l)=1$ and $ \vp(0,\l)=\vt'(0,\l)=0. $
 Here and in the following $u'$ denotes the derivative w.r.t. the first
variable. The Lyapunov function (Hill discriminant) of the periodic
equation \er{1} is then defined by
$$
\gF(\l)={1\/2}(\vp'(1,\l)+\vt(1,\l)),\qqq \l\in \C.
$$
Note that $\s_{ac}(h)=\{\l\in \R: \gF(\l)\in [-1,1]\}$, see
\cite{Ti58}.

\subsection{Riemann surface}
 The function $\gF(\l)$ is entire and is real on the real line.
 Introduce the function $\f$, which is analytic in $\C_+$ and given
 by
\[
\lb{dF}
 \f(\l)=(1-\gF^2(\l))^{1\/2},\qqq \l\in
 \ol\C_+,
\]
 where the branch is defined by the condition
$\f(\l+i0)>0$ for $\l\in \s_0=[\l^+_{0},\l^-_1]$. We also introduce
the two-sheeted Riemann surface $\L$ of $\f(\l)$ obtained by joining
the upper and lower rims of two copies of the cut plane
$\C\sm\s_{ac}(h)$ in the usual (crosswise) way, see e.g.
\cite{KS12}. We denote the $n$-th open gap on the first, physical
sheet $\L_1$ by $\g_n^{(1)}$ and its counterpart on the second,
nonphysical sheet $\L_2$ by $\g_n^{(2)}$, and set the "circle" gap
$\g_n^\bu$ by
\[
\lb{sL}
 \g_n^\bu:=\ol\g_n^{(1)}\cup \ol\g_n^{(2)}.
\]

\subsection{Floquet solutions}
The Floquet solutions $\p^{\pm}(x,\l), \l \in \L$, of the equation
\er{1} are given by
\[
\lb{3}
\begin{aligned}
\p_\pm(x,\l)=\vt(x,\l)+m_\pm(\l)\vp(x,\l),\qqq
\\
\where \qq  m_\pm={a\pm i\f\/ \vp(1,\cdot)},\qqq
a(\l)={1\/2}(\vp'(1,\l)-\vt(1,\l)),
\end{aligned}
\]

$$
\vp(1,\l)\p_+(\cdot,\l)\in L^2(\R_+)\qqq  \forall \ \l\in\L_1=\C\sm
\s_{ac}(h).
$$
Note that in the trivial case $v=0$, we have $\p^\pm(x,\l)=e^{\pm
ix\sqrt \l}$.

Introduce a function $b(\l)=-i\f(\l)$. It is known that (see e.g.
\cite{K05})
\[
\lb{bg}
  b(\l)=(-1)^n\sqrt{\gF^2(\l
)-1}, \ \ \ \ \ \l\in \g _n \in \L ,
 \]
where the branch $\sqrt{\gF(\l)^2-1}>0$ as $\l\in \g _n^1\ss \L^1$
and $\sqrt{\gF(\l)^2-1}<0$ as $\l\in \g _n^2\ss \L^2$.

Below we need the simple identities
\[
\lb{LD0} a^2+1-\gF^2=1-\vp'(1,\cdot)\vt(1,\cdot)=
-\vp(1,\cdot)\vt'(1,\cdot).
\]
%This yields $ m_+(\l)m_-(\l)= -{\vt'(1,\l)\/\vp(1,\l)},\  \l\ne
%\m_n$.

\subsection{Eigenvalues and resonances}
It is well known (see e.g. \cite{KS12}, \cite{Z69})  that, for each
$f\in C_0^\iy(\R_+), f\ne 0$, the function $g(\l)=((h-\l)^{-1}f,f)$
has a meromorphic extension from the physical sheet $\L_1$ to the
whole Riemann surface $\L$. By definition,

If $g$ has a pole at some $\l_0\in \L_1$ for some $f$, then $\l_0$
is an eigenvalue  of $h$ and $\l_0\in\bigcup\limits_{n\ge 0}
\g_n^{(1)}$.

If $g$ has a pole at some $\l_0\in \L_2$ for some $f$, then $\l_0$
is a resonance of $h$. In particular, if $\l_0\in \g_n^{(2)}$ for
some $n\ge 0$, then $\l_0$ is anti-bound state of $h$ and
$\l_0\in\bigcup\limits_{n\ge 1} \g_n^{(2)}$.

If $g$ has asymptotics $g(\l)={1\/\sqrt {t}}(1+O(t))$ as
$t=\l-\l_0\to 0$ for some $\l_0\in\{\l_n^+, \l_n^-\}, n\ge 1$, then
$\l_0$ is a virtual state  of $h$.

It is well known that for the case $v\ne \const$, see \cite{KS12},
\cite{Z69},  the function $g$  has exactly one simple pole
$\m_n^\bu$ on each "circle" gap $\g_n^\bu\ne \es, n\ge 1$ and there
are no others. Here $\m_n^\bu\in \g_n^\bu$ is a so-called state of
$h$ and   its projection onto the complex plane coincides with the
$n$-th eigenvalue, $\m_n$, of  the Dirichlet boundary value problem
$$
-y_n''+vy_n=\m_n y_n, \qq y_n(0)=y_n(1)=0, \qqq x\in [0, 1],\qq n\ge
1.
$$
Moreover, if $\g_n\ne \es$, then  exactly one of the following three
cases holds,

1) $\m_n^\bu\in \g_n^{(1)}$ is an eigenvalue,

2) $\m_n^\bu\in \g_n^{(2)}$ is a resonance (it is a so-called
anti-bound state),

3) $\m_n^\bu\in\{\l_n^+, \l_n^-\}$ is a virtual state. Here the
function $g(\m_n^\bu +z^2)$ has a pole at $0.$

\bigskip

\no There are no other states of $h$, so $h$ has only eigenvalues,
virtual states and anti-bound states. If there are exactly $N\ge 1$
non-degenerate gaps in the spectrum of $\s_{ac}(h)$, then the
operator $h$ has exactly $N$ states; the closed gaps $\g_n=\es$ and
the semi infinite gap $(-\iy,0)$ do not contribute any states. In
particular, if $\g_n=\es$ for all $n\ge 1$, then $v=\const$
(\cite{K98}, \cite{K99}) and thus $h$ has no states.
%The states $\l_n^0$ are described in Lemma \ref{Tm}.
A more detailed description of the states of $h$ is given in Theorem
 \ref{Ti} below.

\subsection{Inverse problem}

We need the following results from the inverse spectral theory for
the operator $h$ on the half-line, in the form convenient for us. We
define the real Hilbert spaces
$$
\mH_\a=\rt\{q\in L^2(0,1): \int_0^1q(x)dx=0, \qq q^{(\a)}\in
L^2(0,1)\rt\}, \a\ge 0,
$$
and let $\mH=\mH_0$.  Introduce  the real Hilbert spaces $\ell_\a^2,
\a\in \R$ of the sequences $(f_n)_1^\iy$ equipped with the norms
$$
\| f\| _{\ell_\a^2}^2=\sum_{n\geq 1}(2\pi n)^{2\a}f_n^2<\iy,
$$
and let $\ell ^2=\ell _0^2$.

Defining the mapping $\x: \mH\to \ell^2\os \ell^2$ by
$$
v\mapsto \x=(\x_n)_1^\iy, \qqq \x_n=(\x_{1n},\x_{2n})\in \R^2
$$
where the components $\x_{1n},\x_{2n}$ are given by
\[
\lb{di}
 \x_{1n}={\l_n^-+\l_n^+\/2}-\m_n^2,\qqq
\x_{2n}=\rt|{|\g_n|^2\/4}-\x_{1n}^2\rt|^{1\/2}\sign_n,\qq
\x_n^2={1\/4}|\g_n|^2,
\]
where
\[
\lb{di1}
 \sign_n=\ca +1  & {\rm if} \ \m_n^\bu \ {\rm is \ an \ eigenvalue,} \\
        -1 & {\rm if} \ \m_n^\bu \ {\rm is \ a \ resonance,}      \\
         0 & {\rm if}  \ \m_n^\bu \ {\rm is \ a \ virtual \ state,}
         \ac,\qqq n\ge 1.
\]
 This mapping $\x$ is described by the following result from
\cite{K99}, \cite{K98}, \cite{K06}:

\begin{theorem}
\lb{Ti} The mapping $v\to \x$ acting from $\cH$ to $\ell^2\os
\ell^2$ is a real analytic isomorphism between the real Hilbert
spaces $\mH$ and $\ell^2\os \ell^2$ and satisfies
\[
\lb{esg} \|v\|\le 4\|\x\| (1+ \|\x\|^{1\/3}),\qqq \|\x\|\le
\|v\|(1+\|v\|)^{1\/3},
\]
where $\|v\|^2=\int_0^1v^2(x)dx$ and $\|\x\|^2={1\/4}\sum_{n\ge 1}
|\g_n|^2$. Moreover, if $\a\in \N$, then $(|\g_n|)_1^\iy\in
\ell_\a^2$ iff the derivative $q^{(\a)}$ belongs to $\mH$.

\end{theorem}

\no{\bf Remark.} 1) Thus we have that for any non-negative sequence
$\vk=(\vk_n)_1^\iy\in \ell^2$, there are unique 2-periodic
eigenvalues $\l_n^\pm$ $(n\in\N_0)$, for some $p\in \cH$, such that
each $\vk_n=\l_n^+-\l_n^-,$ $(n\in\N)$. Consequently, from the gap
lengths $(|\g_n|)_1^\iy$ one can uniquely recover the Riemann
surface $\L$ as well as the points $\l_n^-=\l_n^+$ where $\vk_n=0$.
Furthermore, for any additional sequence $\wt\l_n^0\in \g_n^\bu$
$(n\in\N)$, there is a unique potential $p\in \cH$ such that each
state $\l_n^0$ of the corresponding operator coincides with
$\wt\l_n^0$ $(n\in\N)$. The results of \cite{K99} were extended in
\cite{K03} to periodic distributions $v=p'$, where $p\in\cH$.

2) We can consider the Schr\"odinger operator $\wt h$ acting on the
space $L^2(\R_+ )$ and  given by
$$
hy=-y''+vy, \qqq
$$
with the Neumann boundary condition $y(0)=0$.
 The spectrum of $\wt h$ consists of an
absolutely continuous part $\s_{ac}(\wt h)=\s_{ac}(h)$ (the union of
the bands $\s_n, n\ge 0$ separated by gaps $\g_n$, see \er{sh}) plus
at most one eigenvalue in each non-empty gap $\g_n$, $n\in\N$. Here
also we can consider the resonances for the Schr\"odinger operator
$\wt h$, similar to the the resonances for the Schr\"odinger
operator $h$. These eigenvalues and the resonances for $\wt h$
coincide with the eigenvalues $\n_n, n\ge 0$ of the problem
\[
\lb{Np} -y''+vy=\l y,\qqq y'(0)=y'(1)=0.
\]
The eigenvalues $\n_n, n\ge 0$ satisfy
\[
\lb{Np1} \n_0\le \l_0^+,\qqq \n_n\in [\l_n^-, \l_n^+]\qq \forall \
n\ge 1.
\]
Here also we can consider the resonances for the Schr\"odinger
operator $\wt h$, similar to the the resonances for the
Schr\"odinger operator $h$.

\section {One dimensional  half-solid }
\setcounter{equation}{0}

In this section we consider the case of one-dimensional octant
periodic potentials in the specific form given by \er{qs}.
%We recall the result from \cite{K05}.
We consider the half-solid operator $T_\t, \t\in \R$ acting on
$L^2(\R )$ and given by
\[
\lb{qs}
 \ T_\t=-{d^2\/dx^2}+q_\t(x),\qqq q_\t(x)=\ca \t & if \ x<0\\
           v(x) & if \ x>0\ac,\
\]
where the potential $v\in L_{real}^2(\T)$. By the physical point of
view $v$ is the potential of a crystal and the constant $\t$ is the
potential of a vacuum. In order to describe the spectrum of $T_\t$
we use some properties of the operator $h=-{d^2\/dx^2}+v$ on the
half-line from Section 2.
%Recall that operator $h$ has the spectrum
%given by \er{sh}-\er{00x}.
% We assume that $\l_0^+=0$????.

We recall needed results about  operators $T_\t$ from \cite{K05}.
 We have the following simple results about the spectrum of  $\s (T_\t)$
 given by
\[
\lb{Tac}
\begin{aligned}
\s(T_\t)=\s_{ac}(T_\t)\cup \s_{disc}(T_\t),\qqq
 \s_{ac}(T_\t)=\s_{ac}(h)\cup [\t,\iy).
%& \s_{disc}(T_\t)\ss \R\sm \s_{ac}(T_\t)=(\l_N^-, )\cup
%\bigcup_{j=0,...,N-2} \g_j(h)  ?????
\end{aligned}
\]
Our goal is to study the eigenvalues in the gaps $\g_n(T_\t), n\ge
0,$ and to find how these eigenvalues depend on $\t$. We take any
integer $N\ge 1$. We describe the basic properties of the one
dimensional  half-solid operator $T_\t$:

$\bu $ If $\t\le \l_1^-$, then there is no any gap in the spectrum
of $T_\t$ and we get
\[
\s_{ac}(T_\t)= (\t_0,\iy ), \qqq where \qqq \t_0=\min \{\l_0^+,\t\}.
\]

 $\bu $ If $\t\in \g_N\cup \s_{N}=(\l_N^-, \l_{N+1}^-]$, for some  $N\ge 1$,
 then  the spectrum of $T_\t$ has the form:
\[
\s (T_\t)=\s_0\cup\s_1\cup.... \s_{N-1}\cup\wt\s_N, \qqq
\wt\s_N=[\t_N,\iy),\qq \t_N=\min \{\l_N^+,\t\}.
%\ca [s,\iy)\ & \ if \ s\in (\l_N^-, \l_N^+)\\ [\l_{N}^+,\iy)\ &\ if
%\ s\in \s_{N}\ac.
\]
Thus, there are possible gaps in the spectrum $\s_{ac}(T_\t)$ given
by
\[
\lb{g1}
\begin{aligned}
\g_j(T_\t)=\g_j(h),\qq j=0,1,.., N-1,\qq \g_N(T_\t)=(\l_N^-,\t_N).
%  \s_{disc}(T_\t)\ss \R\sm
%\s_{ac}(T_\t)=(\l_N^-, )\cup \bigcup_{j=0,...,N-2} \g_j(h)  ?????
\end{aligned}
\]
Moreover, in each open gap $\g_j(T_\t)\ne \es, j=0,1,...,N$ there is
at most one eigenvalue $\m_j(\t)$.

\
 $\bu $
 We introduce the Weyl-type functions $\P_{\pm}$,  which are solutions of the equation
\[
\lb{ehs} -y''+q_\t y=\l y
\]
and satisfy
$$
\P_{\pm}(\cdot,\l)\in L^2({\R_{\pm}}), \qqq \forall \ \l\in \L_\t=\C
\sm \s_{ac}(T_\t).
$$
They  have the forms
$$
\begin{aligned}
\P_{+}(x,\l)=\p^+(x,\l), \ \  x \ge 0; \\
\P_{-}(x,\l)=e^{x\sqrt{\t-\l}},  \ \  x \le 0,
\end{aligned}
$$
for $\l<\t, \l\in \L_\t$ and here $\p^+$ given by \er{3}. These
functions $\P_{\pm}(x,\l)$ are analytic in the cut domain $\L_\t$
and are continuous up to the boundary.

$\bu$  We define the Wronskians
\[
w(\l)=\{\P_-,\P_+\}=m_+(\l)-\sqrt{\t-\l}, \ \ \l\in \L_\t.
\]
The function $w(\l)$ on the first sheet $\L_\t$ has finite number of
zeros, which are simple and coincide with eigenvalues of the
operator $T_\t$.

$\bu$
%It is clear that they depend on $t$ periodically.
Using \er{3} we rewrite the Wronskian $w(\l)$ in the gap $\g_n\ss
\L_\t, \l<\t$ in the form
\[
\lb{w1}
\begin{aligned}
w(\l)=m^+(\l)-\sqrt{\t-\l}={a(\l)-b(\l)\/\vp(1,\l)}-\sqrt{\t-\l}, \
\ \ \l<\t,\  \l\in\g_n\ss \L_\t,
\end{aligned}
\]
where
\[
\lb{w2}
\begin{aligned}
\sqrt{\t-\l}>0, \qq {\rm if}\  \ \l<\t, \l\in \L_\t,\qqq
 (-1)^nb(\l)=\sqrt{\D^2(\l
)-1}>0, \qq \l\in \g _n \in \L_\t,
\end{aligned}
\]

\begin{lemma}
\lb{TD1} Assume that the operator $hy=-y''+vy, y(0)=0$ on
$L^2(\R_+)$ has an open gap $I=(\l^-,\l^+)$ in the continuous
spectrum and an eigenvalue $\m\in (\l^-,\l^+)$ for some $v\in
L^2(\T)$. Then for any constant $\t$ large enough the operator
$T_\t$ defined by \er{qs} has an eigenvalue $\m_\t\in I$ such that
\[
\lb{ms} \m_\t-\m={c(\m)\/\sqrt \t}+{O(1)\/\t}  \qqq \as \qqq \t\to
\iy,
\]
where $c(\m)={2b(\m)\/\vp_\l(1,\m)}\ne 0$ and here
$\vp_\l(1,\m)={\pa \/\pa \l}\vp(1,\m)$.
\end{lemma}

\no {\bf Proof.} i) Due to \er{w1} the eigenvalues of $T_\t$ are
zeros of  the Wronskian
$$
w(\l)=m_+(\l)-\sqrt{\t-\l}={a(\l)-b(\l)\/\vp(1,\l)}-\sqrt{\t-\l}
$$ on
the first sheet $\L_\t$. Consider the two functions
$m_+(\l)={a(\l)-b(\l)\/\vp(1,\l)}$ and $\sqrt{\t-\l}$ on the gap
$(\l^-,\l^+)$, where $\t>>\l_+$. The point $\m\in I$ is an
eigenvalue of the operator $h$. Then due to \er{LD0} we have
$a^2(\m)=b^2(\m)\ne 0$. Then the function $m_+(\l)$ is a meromorphic
in the disk around the centrum of the gap $\g_n$ and has the
following asymptotics
\[
m_+(\l)={c(\m)\/\l-\m}+O(1)\qqq \as \qq \l\to \m.
\]
 Thus  the equation
$m_+(\l)=\sqrt{\t-\l}$ has a unique solution $\m_\t\to \m$ as $\t\to
\iy$ given by \er{ms}, since
$$
{c(\m)\/\m_\t-\m}+O(1)=\sqrt \t+O(\t^{-{1\/2}}).
$$
 \BBox

Now we prove the main result of this section.

\begin{lemma}
\lb{TD2} i) Let integer $N\ge 1,\a\ge 0$ and let $\g>0$. Then there
exists a potential $v\in \mH_\a$ such that the first $N$ gaps  in
the spectrum of the operator $h$ on $L^2(\R_+)$ are open and satisfy
\[
|\g_j|=\g,\qq \forall \ j=1,2,...,N.
\]
Moreover, in addition for any points $\l_j\in\g_j, j=1,2,...,N $,
then exists a periodic potential $v\in \mH_\a$ such that each
$\l_j=\m_j, j=1,2,...,N $ is an eigenvalue of the operator $h$.

ii) Let in addition $q_\t$ be given by \er{qs} and let $\t$ be large
enough. Then  each $\g_j, j=1,2,...,N$ is a gap in the spectrum of
$T_\t$ and on each $\g_j$ exists an  eigenvalue $\m_j(\t)\in \g_j$
such that
%$\m_j(s)\to \m_j$ as $s\to \iy$
\[
\lb{mns} \m_j(\t)-\m_j={c(\m_j)\/\t^{1\/2}}+{O(1)\/\t}  \qqq \as
\qqq \t\to \iy,
\]
where $c(\m)={2b(\m)\/\vp_\l(1,\m)}\ne 0$.

iii) Let $\l_0^+=0$ and let $\n_0\le \l_0^+=0$ be the first Neumann
eigenvalue (see \er{Np1}) of the problem \er{Np}. Then for any
$\a\ge 0$ there exists a potential $v\in \mH_\a$ such that $
m_+(0)>0$ and the operator $T_\t$ has an eigenvalue $E<0$ for each
$\t\in (\n_0, m_+(0)^2)\ne \es$.
\end{lemma}

\no {\bf Proof.} i) It follows from  Theorem \ref{Ti}.
%   was proved in \cite{K06} for $\a=0$, see
%Theorem \ref{Ti}.

ii) It follows from Lemma \ref{TD1}.

iii) We define $\r:=m_+(\l_0^+)={a(\l_0^+) \/\vp(1,\l_0^+) }$ for
some potential $v\in \mH_\a$. We recall a needed result from
\cite{K05} about the first eigenvalue:

if $\r<0 $, then $\#( T_\t, \g_0(T_\t))=0 $,

if $\r>0$, then
\[
\lb{w3} \#(T_\t,\g_0(T_\t))=\ca 0, \ \ &if \ \  \t\le \n_0\ {\rm
                                         or}\ \t\ge \r^2 \\
                             1, \ \ &if \ \ \nu_0<\t< \r^2\ac .
\]

Assume that $\r>0$ for some potential $v\in \mH_\a$. Then due to
\er{w3} the operator $T_\t$ has an eigenvalue $E<\l_0^+=0$ for each
$\t\in (\n_0, \r^2)$.

 We show that $\r>0$ for some
potential $v\in \mH_\a$. Below we take $v=p(x+t)$ for some $p\in
L^2(\T)$ and small $t$. We assume that $p$ satisfies\\
{\bf Condition P}. {\it 1) the function $p, p''\in L^2(\T)$;\\
2) $ p(1-x)=p(x), \ \forall \ x\in [0,1]$, i.e., the potential $p$ is even on the interval $[0,1]$ ;\\
3) $p(x)>\d>0$ for all $x\in [0,\ve]$ for some small constants $\d,
\ve>0$ }.

Recall that we put $\l_0^+=0$ and in this case we have
$\int_0^1p(x)dx>0$ (see e.g. \cite{K97}), and then  item 3) in
Condition P is possible since $\int_0^1p(x)dx>0$. We define the
fundamental solutions $\vp(x,\l,t), \vt(x,\l,t)$  of the following
equation with the shifted potential
\[
\lb{eqt}
\begin{aligned}
-y''+p(x+t)y=\l y, \ \ \ \l\in \C ,\ \ \ \  t\in \R ,
\\
\vp_x(0,\l,t)=\vt (0,\l,t)=1,\qqq \vp(0,\l,t)=\vt_x(0,\l,t)=0.
\end{aligned}
\]
 For the shifted potential
$v=p(\cdot+t)$ we define the Lyapunov functions
$\gF(\l,t)={1\/2}(\vp'(1,\l,t)+\vt(1,\l,t))$. Note that we have
$\gF(\l,t)=\gF(\l,0)$, i.e.,  the Lyapunov function $\gF(\l,t)$ for
\er{eqt} coincides with the Lyapunov function for the case $t=0$
(see \cite{L87}). We also define the functions
$$
\begin{aligned}
 a(\l,t)={1\/2}(\vp'(1,\l,t)-\vt(1,\l,t)),\qqq
m_+(\l,t)={a(\l,t)-\f(\l)\/\vp(1,\l,t)}.
\end{aligned}
$$
 Let $\dot u={\pa\/\pa t}u$. We have the equations
\[
\begin{aligned}
\dot a(\l,t)=-\vt_x(1,\l,t)-(\l-p(t))\vp(1,\l,t), \qqq \forall \
(\l,t)\in \C\ts\R,
\end{aligned}
\]
see \cite{L87}.  Then the properties of $p$ give
$$
\dot a(0,t)=-\vt'(1,0,t)+p(t)\vp(1,0,t)>\d\vp(1,0,t)>0\qqq \forall \
t\in [0,\ve],
$$
since $\vt'(1,0,t)\le 0$  (its first zero $\n_0\le 0$ ) and
$\vp(1,0,t)>0$  (its first zero $\m_1>0$) which yields
\[
a(0,t)=a(0,0)+\int_0^t\dot a(0,\t)d\t=\int_0^t\dot a(0,\t)d\t>0,
\]
since $a(0,0)=0$ for all even potentials (see Lemma 3.4. in
\cite{K05}). This implies $m_+(0,t)={a(0,t)\/\vp(1,0,t)}>0$ for the
potential $v(x)=p(x+t)$ all $t\in [0,\ve]$. \BBox

\section {Model operators on $\R_+^d$ and  $\R^d$}
\setcounter{equation}{0}

\subsection{Specific periodic Schr\"odinger operators on the half-line}
 Consider the Schr\"odinger  operator giveb by
 $$
 hf=-f''+vf\qq  {\rm on} \qq L^2(\R_+),\qqq f(0)=0.
 $$
  Recall that the
spectrum of $h$ consists of an absolutely continuous part  (which is
a union of non-degenerate spectral bands $\s_n=[\l^+_n,\l^-_{n+1}],
n\ge 0$) plus at most one eigenvalue in each open  gap
$\g_n=(\l^-_{n},\l^+_n), n\ge 1$ between bands
\cite{KS12}, \cite{Z69} (see Fig.1) and the $\l_n^\pm$ satisfy \er{00x}.

{\it Now we begin to construct a specific potential $v$.} Here we
use results about the gap-lengths mapping from Lemma \ref{TD2} i).
Due to these results about the gap-lengths mapping,
 we take the potential $v\in \mH_\a$ for any fixed $\a\ge 0$ such that the first $N$ gaps
 $\g_1, ..., \g_N$ and other ones $\g_n, n>N$ in the spectrum of $h$    satisfy
\[
\lb{bg1}
\begin{aligned}
&  \g=|\g_1|=|\g_2|=|\g_3|=....|\g_N|= {\pi^24 (N+1)^2\/\vk},\qqq
0<\vk <<1,
\\
&   \sum_{n>N}n^{2\a}|\g_n|^2=Q<\iy.
\end{aligned}
\]
The value $Q$ is not important in our consideration and thus we can
take any $\a\ge 0$. Thus \er{bg1} and the estimate \er{eb} give
\[
\lb{bg1A}
\begin{aligned}
\l_n^-=\g(n-1)+A_n, \qqq \l_n^+=\g n+A_n, \qqq
A_n=\sum_{j=0}^{n}|\s_j|,
\\
|A_{n}|\le \pi^2\sum_{j=0}^{n}(2j+1)=\pi^2(n+1)^2\le \pi^2(N+1)^2\le
{\vk\/4}\g.
\end{aligned}
\]
 Due to Theorem
\ref{Ti} in each open gap $\g_n$, $n=1,2,...,N$ we choose
exactly one eigenvalue $\m_n^\bu$ by
\[
\lb{bg2} \m_n^\bu=\l_n^-+{\g\/4d} \in \g_n^1.
\]
Moreover, \er{bg1A} gives
\[
\lb{bg3}
\begin{aligned}
\textstyle
 \m_n^\bu=\l_n^-+{\g\/4d}=\g(n-1+{1\/4d})+A_n,
\qqq \forall \ n=1,..,N.
\end{aligned}
\]
It is convenient to define the "normalized" operator
$$
h_\g={1\/\g}h.
$$
Then the spectrum of $h_\g$ consists of an absolutely continuous
part $\s_{ac}(h_\g)= \bigcup\limits_{n\ge 0} s_n $ plus at most one
eigenvalue in each non-empty gap $g_n$, $n\in\N$, but exactly one
eigenvalue $e_n={\m_n\/\g}$ in each open gap $g_n, n=1,...,N$. Here
the spectrum of $h_\g$ has  the bands $s_n={\s_n\/\g}$ and gaps
$g_n={\g_n\/\g}$. In particular, we have
$$
\textstyle s_0={\s_0\/\g},\qqq
s_n={\s_n\/\g}=[n+{A_n\/\g},n+{A_{n+1}\/\g}],\ \ \qq
g_{n}={\g_n\/\g}, \qq n=1,2,...,N,
$$
where $A_n$ is defined in \er{bg1A}. The first bands  satisfy
\[
\lb{ebx} \textstyle
 |s_n|={|\s_{n}|\/\g}\le {\pi^2\/\g}(2n+1)\le {\vk\/2(N+1)}
%:={\pi^2\/\g}(2N+1)={1 \/\a}
,\qqq \forall \  \qqq n=0,...,N.
\]
Thus these spectral bands $s_n$ are very small and very close to the
points $n$ and satisfy
\[
\lb{snn} \dist \{s_n,n\}\le {A_{n+1}\/\g}\le {\vk\/4}, \qq
n=0,1,...,N.
\]
In this case due to \er{bg1} the gaps $g_n$ satisfy
\[
\lb{es1}
\begin{aligned}
|g_n|=1, \qqq n=1,...,N.
\end{aligned}
\]
In each gap $g_n, n=1,...,N$, there exists exactly one eigenvalue
$e_n$ of $h_\g$ such that
\[
\lb{es2}
\begin{aligned}
\textstyle
 e_n={\m_n\/\g}=e_n^0+{A_n\/\g},\qqq
e_n^0=n-1+{1\/4d}=n-1+e_1,\qq e_1={1\/4d},\qqq {A_n\/\g}\le
{\vk\/4}.
\end{aligned}
\]
%for all  $n=1,...,N$.

%\newpage

\subsection{\lb{H2} Schr\"odinger
operators on $\R_+^2$}
  We consider  Schr\"odinger operators
  $H_+=h_1+h_2$ on the quadrant  $\R_+^2$.
Here $h_1$ and $h_2$ are defined on the half-line and  depend on one
variable and given by
$$
h_jy=-y''+v(x_j)y, \qqq y(0)=0, \qq j=1,2,
$$
where the potential $v\in L^2(\T)$. For large constant $\g$ we
define a new operator
$$
H_\g={1\/\g}H_+={1\/\g}(h_1+h_2).
$$
We take the operator $H_\g$, when the variables are separated.  We
show that $ H_\g$ has (first) bands which are very small and their
positions are very close to the integer $n$. The union of group of
bands close to the integer $n$ is a cluster $K_n$. Between the two
neighbor clusters $K_n$ and $K_{n+1}$ there exists a big gap. On
this gap there exist $n$ eigenvalues.

$\bu$ We define bands (i.e., the basic bands) $S_{i,j}^0$ of the
operator ${H_\g}$ and their clusters $K_n^0$ by
\[
S_{i,j}^0=s_i+s_j,\qqq \qqq K_n^0=\bigcup_{i+j=n} S_{i,j}^0,\qqq
i,j,n\in \Z_+=\{0,1,2,3,....\},
\]
where we define $A+B=\{z=x+y:x\in A, b\in B\}$ for sets $A,B$.  In
particular, we have
\[
\begin{aligned}
K_0^0=S_{0,0}^0,\qq K_1^0=S_{1,0}^0,\qq K_2^0=S_{2,0}^0\cup
S_{1,1}^0,.....,
%\\K_n=\bigcup_{i+j=n} (s_i+s_j),....,\\
\end{aligned}
\]
If $\g$ is large enough, then due to \er{snn}, \er{ebx} we estimate
the position of bands $S_{i,j}^0$ and their lengths $|S_{i,j}^0|$ by
\[
\lb{eS0} \dist \{S_{i,j}^0, i+j\}\le {\vk\/2},\qqq |S_{i,j}^0|\le
{\vk\/(N+1)}
\]
for all $i,j=0,1,..., N$. This yields the position of their clusters
$K_n^0$ and their diameters by
\[
\lb{eK0} \dist \{K_n^0, n\}\le {\vk\/2} \qqq \diam K_n^0\le
{3\/2}\vk
\]
for all $n=0,1,..., N$.

%  These clusters  are
%separated by gaps $G_n=(K_{n-1}^+, K_n^-), n\ge 1$. Thus we have

%From \er{} we deduce that the following
%\[
%K_n\ss [K_n^-,K_n^+], \qq K_n^-=\min_{i+j=n} (\l_i^-+\l_j^-), \qq
%K_n^+=\max_{i+j=n} (\l_i^++\l_j^+), G_n\ss , n=1,2,...,N
%\]

$\bu$ We define the surface bands $S_{i,j}^1$ of the operator  ${H_\g}$
and their clusters $K_n^1, $ by
\[
\lb{dS1} S_{i,j}^1=e_i+s_j,\qqq \qqq K_n^1=\bigcup_{i+j=n+1}
S_{i,j}^1, \qq i,j, n\ge 0.
\]
In particular, we have
%These clusters  are separated by gaps $\cG_n,
%n=1,2,...,N$.  Thus we have
\[
\lb{dK1}
\begin{aligned}
K_0^1=S_{0,1},\qq K_1^1=S_{0,2}^1\cup S_{1,1}^1,\qq
 K_1^1=S_{0,3}^1\cup S_{1,2}^1\cup S_{2,1}^1,....
%\mS_N=S_{0,N}^1\cup S_{1,N-1}^1\cup.....
\end{aligned}
\]
Using arguments similar to the case of the bands $S_{i,j}^0$ we
determine the position of surface  bands $S_{i,j}^1$ and their
lengths $|S_{i,j}^1|$ by
\[
\lb{eE1} \dist \{S_{i,j}^1, e_i^0+j\}\le {\vk\/2},\qqq
|S_{i,j}^1|\le {\vk\/2(N+1)}
\]
for all $i,j=0,1,..., N$. This yields the position of their clusters
$K_n^1$ and their diameters by
\[
\lb{eK1} \dist \{K_n^1,e_1^0+n\}\le {\vk\/2}, \qqq \diam K_n^1 \le
\vk
\]
for all $n=0,1,..., N$.

$\bu$ The operator $H_\g$ has {\bf eigenvalues} $E_{i,j} $ and their
cluster $K_n^e$  given by
\[
E_{i,j}=e_i+e_j,\qqq i,j\ge 1,\qqq K_n^e=\{E=E_{i,j}, i+j=n+1\},\qq
n\ge 1,
\]
The cluster $K_n^e$  has n eigenvalues of the operator $H_+$. In
particular, we have
\[
\begin{aligned}
\lb{Ke1} K_1^e=\{E_{1,1} \},\qqq K_2^e=\{E_{1,2},\ E_{2,1} \},\qqq
K_3^e=\{E_{1,3},\ E_{2,2}, \ E_{3,1}\},....
\end{aligned}
\]
The identity \er{es2} gives
\[
\begin{aligned}
\lb{Ke2} \textstyle E_{i,j}=e_i^0+e_j^0+{A_i\/\g}+{A_j\/\g},\qqq
|E_{i,j}-e_i^0-e_j^0|\le {\vk\/2}
\end{aligned}
\]
for all $i,j=1,2,...,N$. This yields the position of their clusters
$K_n^e$ and their diameters by
\[
\lb{Ke3} \textstyle \dist \{K_n^e,2e_1+ n-1\}\le {\vk\/2}, \qqq
\diam K_n^e\le \vk,
\]
for all $n=1,2,..., N$.

$\bu$ Thus we can describe  $\s_{ac}(H_+)$ and $\s_{disc}(H_+)$ by
\[
\s_{ac}(H_+)=\cup_{n\ge 0} (K_n^0\cup K_n^1),\qqq
\s_{disc}(H_+)=\cup_{n\ge 1}K_n^e.
\]
Now  combining all estimates \er{eK0}-\er{eK1} we deduce that
between two sets $K_n^0\cup K_n^1$ and $K_{n+1}^0\cup K_{n+1}^1$ for
each $n=0,1,...,N$ there exists an interval $I_n$ given by
\[
\lb{Inx}
\begin{aligned}
I_n=[I_n^-, I_n^+]=[e_1^0+n+4\vk, n+1-4\vk],\\
{\rm such \ that}\qqq \dist\{I_n, \s_{ac}(H_+)\}\ge 2\vk,\qqq
\end{aligned}
\]
%where $K=\cup_{n\ge 0} (K_n^0\cup K_n^1)=\s_{ac}(H_+)$,
i.e., the distance between the interval $I_n$ and two sets
$K_n^0\cup K_n^1$ and $K_{n+1}^0\cup K_{n+1}^1$ is greater than
$2\vk$. Moreover, due to \er{Ke1}-\er{Ke3} the eigenvalue cluster
$K_n^e$ satisfies:
\[
\lb{KessI}
\begin{aligned}
K_n^e\ss I_n,\qqq \forall \ n=1,...,N.
\end{aligned}
\]

\subsection{\lb{H3} Schr\"odinger
operators on $\R_+^3$}
We consider  Schr\"odinger operators
$H_+$ on the corner $\R_+^3$ given by
\[
H_+=-\D_++V_+,\qqq V_+(x)=v(x_1)+v(x_2)+v(x_3), \qqq x=(x_j)_1^3\in
\R_+^3
\]
where the potential $v$ is 1-periodic and $v\in L^2(0,1)$. We
rewrite the operator $H$ in the form
\[
H_+=h_1+h_2+h_3, \qq h_j=h_0+v(x_j)
\]
Define the operator $H_\g$ by
$$
H_\g={1\/\g}H=h_{1,\g}+h_{2,\g}+h_{3,\g},\qqq h_{j,\g}={h_j\/\g}.
$$

$\bu$ We define bands $S_{i,j,k}^0$ of the operator  ${H_\g}$ and
their clusters $K_n^0, n=0,1,....,N$ by
\[
S_{i,j,k}^0=s_i+s_j+s_k,\qqq i,j,k\in \Z_+,\qqq
K_n^0=\bigcup_{i+j+k=n} S_{i,j,k}^0,\qqq n\in \Z_+,
\]
and in particular,
$$
\begin{aligned}
& K_0^0=S_{0,0,0}^0=s_0+s_0+s_0,\qqq K_1^0=S_{0,0,1}^0,\qqq
K_2^0=S_{0,0,2}^0\cup S_{0,1,1}^0,....
\end{aligned}
$$
Recall that  we define $A+B$ for sets $A,B$ by $A+B=\{z=x+y:
(x,y)\in A\ts B\} $. Similar to 2dim case  we deduce that
\[
S_{i,j,k}^0\sim i+j+k, \qqq K_n^0 \sim n,\qqq  \forall\qq
n=1,2,...,N.
\]

In 3-dimensional case we have two types of the surface (guided) bands   $S_{i,j,k}^1$ and $S_{i,j,k}^2$.

$\bu$  {\bf The first  type of surface (guided) bands.} We define
the surface (guided) bands $S_{i,j,k}^1$ of the operator ${H_\g}$
and their clusters $K_n^1, n=1,....,p$ by
\[
S_{i,j,k}^1=s_i+s_j+e_k,\qqq \qqq K_n^1=\bigcup_{i+j+k=n+1}
S_{i,j,k}^1,\qqq i,j, n\in \Z_+,\qq k\in \N,
\]
The position of surface bands $S_{i,j,k}^1$ and their clusters
$K_n^1$ are given by
\[
S_{i,j,k}^1\sim i+j+k-1+e_1=n-1+e_1, \qqq \qqq K_n^1\sim n-1+e_1.
\]
These clusters  are separated by gaps $\cG_n, n=1,2,...,p$.  Thus we
have
\[
\begin{aligned}
K_1^1=S_{0,1}^1,\qq K_2^1=S_{0,2}^1,\qq K_3^1=S_{0,3}^1\cup
S_{1,2}^1,...., K_N^1=S_{0,N}^1\cup S_{1,N-1}^1\cup.....
\end{aligned}
\]

$\bu$ {\bf The second type of surface (guided) bands.} We define the
surface (guided) bands $S_{i,j,k}^2$ of the operator ${H_\g}$ and
their clusters $K_n^2, n=1,....,p$ by
\[
S_{i,j,k}^2=e_i+e_j+s_k,\qqq \qqq K_n^2=\bigcup_{i+j+k=n+2}
S_{i,j,k}^2,\qq i,j\in \N,\qq k,n\in \Z_+.
\]
The positions of the surface bands $S_{i,j,k}^2$ and the cluster
$K_n^2$ are given by
\[
S_{i,j}^2\sim i+j+k-2+2e_1=n-2+2e_1, \qqq\qqq K_n^2\sim n-1+2e_1.
\]
These clusters  are separated by gaps $\cG_n, n=1,2,...,p$.  Thus we
have
\[
\begin{aligned}
K_1^2=S_{0,1}^2,\qq K_2=S_{0,2}^2,\qq K_3^2=S_{0,3}^2\cup
S_{1,2}^2,....,\qq K_N=S_{0,N}^2\cup S_{1,N-1}^2\cup.....
\end{aligned}
\]

$\bu$ {\bf Eigenvalues.} The operator $H_\g$ has eigenvalues
$E_{i,j,k} $ and their cluster $K_n^e$ given by
\[
E_{i,j,k}=e_i+e_j+e_k,\qqq i,j,k\in \N,\qqq K_n^e=\{E=E_{i,j,k},
i+j+k=n+3\},\qq n\in\Z_+.
\]
The positions of eigenvalues $E_{i,j,k}, i+j+k=n+3$ and their
cluster $K_n^e$ are given by
\[
E_{i,j,k}\sim i+j+k-3+3e_1=n+3e_1,\qqq K_n^e\sim n+3e_1.
\]
The sets $\s_{ac}(H)$ and $\s_{disc}(H)$ are given by
\[
 \s_{ac}(H_+)=\cup_{n\ge 0} (K_n^0\cup K_n^1\cup K_n^2),\qqq
\s_{disc}(H_+)=\cup_{n\ge 1}K_n^e
\]
Later on we repeat the proof for the case $d=2$. \BBox

\subsection{\lb{h} Specific 1dim half-solid  potentials}
 Consider the  operator $T_\t f=-f''+q_\t f$ on $L^2(\R)$, where the
 potential $q_\t$ is given by \er{qs}.
We take any fix  integer $N\ge 1$ and numbers $\g, \t>0$ large
enough. Due to Lemma \ref{TD2} we obtain that there exists a
periodic potential $v$ such that the first $N$ gaps $\g_j,
j=1,2,...,N$ in the spectrum of the operator $h$ are open. Moreover,
there exists an eigenvalue $\m_j(\t)$ in each this gap $\g_j$ and
they satisfy
\[
\lb{gjxx} \m_j(\t)\in \g_j \qqq |\g_j|=\g,\qq \forall \ j=1,2,...,N,
\]
%The spectrum of $T_s
\[
\lb{Ts1} \s_{ac}(T_\t)=\s_0\cup\s_1\cup.... \s_{N-1}\cup\wt\s_N,
\qqq \wt\s_N\ss [\l_N^+,\iy).
%\ca [s,\iy)\ & \ if \ s\in (\l_N^-, \l_N^+)\\ [\l_{N}^+,\iy)\ &\ if
%\ s\in \s_{N}\ac.
\]
Here the bands $\s_0, \s_1,.... \s_{N-1}$ and the set $\wt\s_N$ are
separated by gaps $\g_j,  j=1,2,...,N$ and each eigenvalue
$\m_n(\t)$  satisfies \er{mns}.

{\bf Now we begin to construct a specific potential $v$.} Here we
use results about the gap-lengths mapping from \cite{K99}. Due to
Theorem \ref{Ti} about the gap-lengths mapping,
 we take the potential $v\in L^2(\T)$ such that the first $N$ gaps
 $\g_1, ..., \g_N$ in the spectrum of $h$ are
 open and   these gaps   satisfy
\[
\g=|\g_1|=...=|\g_N|={\pi^24(N+1)^2\/\vk},\qqq 0<\vk <<1,
\]
and in each big gap $\g_n$, $n=1,2,...,N$ there exists exactly one
eigenvalue $\m_n(\t)$.

 Define the operator
$ T_{\t,\g}={1\/\g}T_\t $. From the properties of $T_\t$ we deduce
that the spectrum of $T_{\t,\g}$ consists of an absolutely
continuous part $\s_{ac}(T_{\t,\g})= \bigcup\limits_{n\ge 0} s_n $
plus  at most one eigenvalue in each non-empty gap $g_n$, $n\in\N$,
 where the bands $s_n$ and gaps $g_n$ are given by
$$
s_0={\s_0\/\g},\qqq  s_n={\s_n\/\g},\qq  g_{n}={\g_n\/\g}, \qq n\in
\N_N=\{1,...,N\},
$$
and they satisfy \er{ebx}-\er{es1}. In each gap $g_n, n=1,...,N$,
there exists exactly one eigenvalue $e_n$ given by
\[
\lb{enx} e_n={\m_n\/\g}=e_n^0+\ve_n,\qqq e_n^0=n-1+{1\/4d},\qq
|\ve_n|\le {\vk\/4},\qqq     n\in \N_N,
\]
since we take $\t$ large enough. Thus roughly speaking the spectrum
of the operators $\T_\t$ on $L^2(\R)$ and $h$ ( on $L^2(\R_+)$) is
the same on the interval $[0,\l_N^+]$. They have the same bands
$\s_0,...,\s_N$ and the same gaps $\g_1,...,\g_N$. Moreover, their
eigenvalues in each gap $\g_n$ are very close, since we take $\t$
large enough.

%\newpage

\subsection{\lb{H2R} Model  Schr\"odinger
operators on $\R^2$}
  We consider  Schr\"odinger operators $H$ on the plane  $\R^2$ given by
\[
  H=T_{\t,1}+T_{\t,2}, \qqq
\]
the proof for the case $\R^d, d\ge 3$ is similar. Here each
$T_{\t,j}=-{d^2\/dx_j^2}+q_{\t}(x_j)$ acts on $\R$ and the potential
$q_{js}$ is determined in Subsection \ref{h}. The spectrum of
$T_{\t,j}$ and $h_j$ are similar on the interval $[0,\l_N^+]$. Then
the spectrum  of the sum $T_{\t,1}+T_{\t,2}$ is similar to the
spectrum of $h_1+h_2$ on the interval $[0,2\l_N^+]$. The proof
repeats the case $h_1+h_2$.

\subsection{ Schr\"odinger
operators on $\R_+\ts \R^{}$} Consider the operator $H=h_1+T_{\t,2}$
on the half-plane $\R_+\ts \R$, where the operator
$h_1y=-y''+v(x_1)y, \ y(0)=0$ acts on the half-line and depends on
one variable $x_1>0$; the operator
$T_{\t,2}=-{d^2\/dx_2^2}+q_{s}(x_2)$ acts on $\R$ and the potential
$q_{\t}(x_2)$ is defined by  \er{qs} and the constant $\t$ is large
enough. The spectrum of $T_{\t,2}$ and $h_1$ are similar on the
interval $[0,\l_N^+]$ for $N, \t$ large enough. The proof repeats
the case $h_1+h_2$.

\section {Proof of main Theorems}
\setcounter{equation}{0}

{\bf Proof Theorem \ref{T1}} i) We consider an operator $H_+=-\D+V$
on $\R_+^2$, where $V$ is $\Z^2$-periodic, the proof for other cases
is similar. Let $H=-\D+V$ on $\R^2$. Define  functions $g_n\in
C_0^\iy(\R)$ and $G_n\in C_0^\iy(\R^2), n\ge 1$ by:
\[
\lb{gn}
\begin{aligned}
g_n|_{w_n}=1, \qq w_n=[4^n, 4^n+n],  \qq \supp g_n=[4^n-1,
4^n+n+1],\\
G_n(x)=g_n(x_1)g_n(x_2),\qqq \supp G_n\ss \R_+^2.
\end{aligned}
\]
Let $\cT_2=\R^2/\Z^2$. For any $\l\in \s(H)$ there exists a function
$\p(x,k)=e^{i(k,x)}u(x,k)$, which  satisfies
\[
\lb{pnu}
\begin{aligned}
(-\D +V(x))\p(x,k)=\l\p(x,k),\qq \forall \ x\in \R^2,\\
u(\cdot,k)\in L^2(\cT_2),\qqq \int_{\cT_2}|u(x,k)|^2dx=1,
\end{aligned}
\]
for some $k\in \R^2$. Define the sequence
$\p_n(x,k)={1\/c_n}G_n(x)\p(x,k)$, where $c_n>0$ is given by
$$
c_n^2=\int_{\R^2}|G_n(x)\p(x,k)|^2dx.
$$
The function $u(x,k)$ is $\Z^2$ periodic, then due to \er{pnu}  we
obtain
\[
c_n^2=\int_{\R^d}|G_n(x)\p(x,k)|^2dx=n^2 +O(n)
\]
as $n\to \iy$. Thus the sequence $\p_n$ satisfies

1) $\|\p_n(\cdot,k)\|=1$ and $\D \p_n\in L^2(\R^2)$, for all $n\in
\N$,

2) $\p_n \perp\p_m$ for  all $n\ne m$, and $\p_n\to 0$ weakly as
$n\to \iy$

Thus $\l\in \s_{ess}(H_+)$, since standard arguments imply
$$
\|(H_+-\l)\p_n(\cdot,k)\|=\|(H-\l)\p_n(\cdot,k)\|\to 0\qq  as \qq
n\to \iy.
$$

ii). Consider an operator $H_\ve=H_++\ve W$ on $\R_+^d$  for the
case $d=2$, the proof for the case $z=x, d\ge 3$ is similar. Here
the operator $H_+$ is defined in subsection \ref{H2}. Recall that
for for any $n\ge1$ there exists a specific potential $v\in L^2(\T)$
such that on the interval $I_n=(a,b)$ (defined by \er{Inx}) contains
$n$ eigenvalues of the operator $H_+$. Moreover, the distance
between the interval $I_n$ and two cluster spectral sets $K_n^0\cup
K_n^1$ and $K_{n+1}^0\cup K_{n+1}^1$ is greater than $2\vk$.

We have $H_0=H_+$, where the real coupling constant $\ve$ is small
enough and  $W$ satisfies
\[
\lb{cW} W\in L^\iy(\R_+^2),\qqq  \|W\|_{L^\iy(\R_+^2)}\le 1, \qqq
0<\ve\le \vk^3.
\]
We define contours $c_n=\{\l\in \C: \dist\{\l,I_n\}=\vk \}$. Due to
\er{cW} the operator $H_\ve$ has $n$ eigenvalues inside the contours
$c_n$, since we have
\[
\begin{aligned}
P_n(\ve)=-{1\/2\pi i}\int_{c_n}R_\ve(z)dz,\\
R_\ve(z)=R_+(z)-R_+(z)\ve W R_\ve(z),\\
\|R_\ve(z)-R_+(z)\|\le \ve \|R_+(z)\| \|R_\ve(z)\|\le
{\ve\/\vk^2}\le \vk\qqq  \forall \ z\in c_n.
\end{aligned}
\]
This yields
\[
\|P_n(\ve)-P_n(0)\|\le \vk<1.
\]
and then the operators  $P_n(\ve)$ and $P_n(0)$ have the same
dimension. In order to show that the intervals $[a_-,a]$ and $[b,
b_+]$ contain the essential spectrum of $H$ for some $a_-<a$ and
$b<b_+$ we use similar arguments.

Now we consider  Schr\"odinger operators $H_\ve=H_0+\ve W$ on the
domain $D=\R_+^{d_1}\ts \R^{d_2},  d_1+d_2=d\ge 2$, where the real
coupling constant $\ve$ is small enough and  $W$ satisfies
\[
\lb{cWz} W\in L^\iy(D),\qqq  \|W\|_{L^\iy(D)}\le 1.
%, \qqq 0<\ve\le\vk^3.
\]
The operator $H_0=H_{01}+H_{02}$, where $H_{01}, H_{02}$ are  given
by
$$
\begin{aligned}
 H_{01}=-\D_x+\sum_{j=1}^{d_1} v(x_j),\qqq  H_{02}=-\D_y+ \sum_{j=1}^{d_2} v(y_j),\qqq
x=(x_j)\in \R_+^{d_1}, y=(y_j)\in \R^{d_2}.
\end{aligned}
$$
The proof for this case is similar and is based on the
 \BBox

{\bf Proof Theorem \ref{T2}.}  We consider  Schr\"odinger operators
$H=T_{\t,1}+T_{\t,2}$ on the plane  $\R^2$, the proof for the case
$\R^d, d\ge 3$ is similar. Here each operator $T_{\t,j}$ acts on
$\R$ and given by
 $$
  T_{\t,j}=-{d^2\/dy_j^2}+q_{\t}(y_j),\qqq y=(y_1,y_2)\in \R^2,
  $$
where the potential $q_{\t}$ is determined in \er{qs}. By Lemma
\ref{TD2}, iii),  each operator $T_{\t,j}, j=1,2$ on $L^2(R)$ has an
eigenvalue $E_0$ below the continuous spectrum
$\s_{ac}(T_{\t,j})=[\t_0,\iy)$ for some potential $v\in L^2(\T)$ and
some $\t\in \R$, where $\t_0=\min\{\t, \l_0^+\}$. Then the operators
$H=T_{\t,1}+T_{\t,2}$ has the eigenvalue $E_0$ below the continuous
spectrum $\s_{ac}(H)=[\t_0+E_0,\iy)$. Consider an operator
$H_\ve=H+\ve W$, where $\ve>0$ is small enough and $W$ satisfies
\er{VO}. It is well know that under small perturbation, the isolated
eigenvalue is still eigenvalue.\BBox

%\newpage

\

\footnotesize \no\textbf{Acknowledgments.} \footnotesize Various
parts of this paper were written during Evgeny Korotyaev's stay as a
VELUX Visiting Professor at the Department of Mathematics, Aarhus
University, Denmark. He is grateful to the institute for the
hospitality. In addition, our study was supported by
%the RSF grant No 15-11-30007 and
the Danish Free Science Council grant
No1323-00360.

\end{document}